\documentclass[12pt]{article}
\usepackage[a4paper, total={7in, 10in}]{geometry}
\usepackage[parfill]{parskip}
\usepackage{physics, tensor, float, subcaption}
\usepackage{graphicx}
\graphicspath{ {Plots/} }
\usepackage{jhep-mod}
\usepackage{bm}
\usepackage{soul}
\usepackage{amssymb,amsmath,amsthm}
\usepackage{mathrsfs}
\usepackage[utf8]{inputenc}
\usepackage{enumerate}
\usepackage{bigints}
\usepackage{xcolor}
\usepackage{appendix}
\usepackage{graphicx}
\usepackage{float}
\usepackage{tikz}
\usepackage{setspace}
\usepackage{cancel}
\usepackage{doi}
\definecolor{purple}{rgb}{1,0,1}
\definecolor{lime}{HTML}{A6CE39} 

\definecolor{lime}{HTML}{A6CE39}
\newcommand{\orcidicon}{%
	\begin{tikzpicture}
	\draw[lime, fill=lime] (0,0) 
		circle [radius=0.16] 
		node[white] {{\fontfamily{qag}\selectfont \tiny ID}};
	\draw[white, fill=white] (-0.0625,0.095) 
		circle [radius=0.007];
	\end{tikzpicture}
	\hspace{-5mm}
}

\newcommand\orcidMatt{{\href{https://orcid.org/0000-0003-1088-6485}{\orcidicon}}}

\renewcommand{\O}{\mathcal{O}}
\begin{document}

\title{\null\vspace{-75pt}\huge{
Effective exponential bounds on the prime gaps\\
}}

\author{
\Large
Matt Visser\!\orcidMatt\!
}
\affiliation{School of Mathematics and Statistics, Victoria University of Wellington, \\
\null\qquad PO Box 600, Wellington 6140, New Zealand.}
\emailAdd{matt.visser@sms.vuw.ac.nz}
\def\theta{\vartheta}
\def\O{{\mathcal{O}}}

\abstract{
\vspace{1em}

Over the last 50 years a large number of effective exponential bounds on the first Chebyshev function $\vartheta(x)$ have been obtained. Specifically we shall be interested in effective exponential bounds of the form
\[
|\vartheta(x)-x| <  a \;x \;(\ln x)^b \; \exp\left(-c\; \sqrt{\ln x}\right); 
\qquad (x \geq x_0).
\]
Herein we shall convert these effective bounds on $\vartheta(x)$ into effective exponential bounds on the prime gaps $g_n = p_{n+1}-p_n$.  Specifically we shall establish
a number of effective exponential bounds of the form
\[
{g_n\over p_n}  < { 2a  \;(\ln p_n)^b \; \exp\left(-c\; \sqrt{\ln p_n}\right) \over
1- a  \;(\ln p_n)^b \; \exp\left(-c\; \sqrt{\ln p_n}\right)};
\qquad (x \geq x_*);
\]
and
\[
{g_n\over p_n}  <  3a  \;(\ln p_n)^b \; \exp\left(-c\; \sqrt{\ln p_n}\right);
\qquad (x \geq x_*);
\]
for some effective computable $x_*$. 
It is the explicit presence of the exponential factor, with known coefficients and known range of validity for the bound,  that makes these bounds particularly interesting.

\bigskip

\bigskip
\noindent
{\sc Date:} 12 November 2022; \LaTeX-ed \today


}

\maketitle
\def\tr{{\mathrm{tr}}}
\def\diag{{\mathrm{diag}}}
\def\cof{{\mathrm{cof}}}
\def\pdet{{\mathrm{pdet}}}
\def\d{{\mathrm{d}}}
\parindent0pt
\parskip7pt
\def\Kerr{{\scriptscriptstyle{\mathrm{Kerr}}}}
\def\eos{{\scriptscriptstyle{\mathrm{eos}}}}
\clearpage
\section{Introduction}

The last 50 years have seen the development of a large number of  fully effective exponential bounds 
on the first Chebyshev function $\vartheta(x)$ --- bounds 
of the form:
\begin{equation}
|\theta(x)-x| <  a \;x \;(\ln x)^b \; \exp\left(-c\; \sqrt{\ln x}\right); 
\qquad (x \geq x_0).
\end{equation}
See references~\cite{Schoenfeld,Trudgian,Johnston-Yang,Fiori-et-al,Broadbent-et-al}.
Here $a>0$ always, while typically $b\geq 0$, and $c>0$ always. 
The special case $b=0$ corresponds to effective bounds of the de la Valle Poussin form~\cite{Poussin,dlVP}. 
For some widely applicable effective bounds of this form see Table I.
(An elementary computation is required for the numerical coefficients in the Schoenfeld~\cite{Schoenfeld} and Trudgian~\cite{Trudgian} bounds.)

\begin{table}[!h]
\caption{Some widely applicable effective bounds 
on the first Chebyshev function $\vartheta(x)$.}\smallskip
\begin{center}
\begin{tabular}{||c|c|c|c|c||}
\hline
\hline
\hline
$a$ & $b$ & $c$ & $x_0$ & Source \\
\hline
\hline
\hline
 0.2196138920& 1/4 & 0.3219796502 & 101 & Schoenfeld~\cite{Schoenfeld}\\
\hline
1/4 & 1/4 & 1/4 & 31 & relaxed version of Schoenfeld~\cite{Schoenfeld}\\
\hline
\hline
0.2428127763 & 1/4 &0.3935970880 & 149 & Trudgian~\cite{Trudgian}\\
\hline
1/4 & 1/4 & 1/3 & 43 & relaxed version of Trudgian~\cite{Trudgian}\\
\hline
\hline
9.220226 & 3/2 & 0.8476836 & 2 & Fiori--Kadiri--Swidinsky~\cite{Fiori-et-al}\\
\hline
\hline
9.40 & 1.515 & 0.8274 & 2 & Johnston--Yang~\cite{Johnston-Yang}\\
\hline
\hline 
\hline
 0.3510691792&0&  1/4  & 101 & 
 Visser~\cite{dlVP}; based on Schoenfeld~\cite{Schoenfeld}\\
\hline
\hline
0.2748124978 &0&  1/4 &  149 & 
Visser~\cite{dlVP}; based on  Trudgian~\cite{Trudgian}\\
\hline
0.4242102935 &0&  1/3 &  149 & 
Visser~\cite{dlVP}; based on Trudgian~\cite{Trudgian}\\
\hline
\hline
295 &0&  1/2 & 2 &  
Visser~\cite{dlVP}; based on FKS~\cite{Fiori-et-al}\\
\hline
\hline
385 &0&  1/2  & 2 & Visser~\cite{dlVP}; based on 
JY~\cite{Johnston-Yang}\\
\hline
\hline
\hline
1 & 0 & 1/4 & 2 & Visser~\cite{dlVP}\\
\hline
1 & 0 & 1/3 & 3 & Visser~\cite{dlVP}\\
\hline
\hline
1/2 & 0 & 1/4 & 29 & Visser~\cite{dlVP}\\
\hline
1/2 & 0 & 1/3 & 41 & Visser~\cite{dlVP}\\
\hline
\hline
\end{tabular}
\end{center}
\end{table}

For some asymptotically more stringent effective bounds, but valid on more restricted regions see Table~II, (based on reference~\cite{Johnston-Yang}), 
and Table~III, (based on reference~\cite{dlVP}).
See also the extensive tabulations in reference~\cite{Broadbent-et-al}.

\clearpage
\begin{table}[!htb]
\caption{Asymptotically stringent bounds on the first Chebyshev function $\vartheta(x)$ valid on restricted regions~\cite{Johnston-Yang}.}\smallskip
\begin{center}
\begin{tabular}{||c|c|c|c||}
\hline
\hline
$a$ & $b$ & $c$ & $x_0$  \\
\hline
\hline
8.87 & 1.514 & 0.8288 & $\exp(3000)$ \\
8.16 & 1.512 & 0.8309 & $\exp(4000)$  \\
7.66 & 1.511 & 0.8324 & $\exp(5000)$  \\
7.23 & 1.510 & 0.8335 & $\exp(6000)$  \\
7.00 & 1.510 & 0.8345 & $\exp(7000)$  \\
6.79 & 1.509 & 0.8353 & $\exp(8000)$  \\
6.59 & 1.509 & 0.8359 & $\exp(9000)$  \\
6.73 & 1.509 & 0.8359 & $\exp(10000)$ \\
\hline\hline
23.14 & 1.503 & 0.8659 & $\exp(10^5)$  \\
38.58 & 1.502 & 1.0318 & $\exp(10^6)$  \\
42.91 & 1.501 & 1.0706 & $\exp(10^7)$  \\
44.42 & 1.501 & 1.0839 & $\exp(10^8)$  \\
44.98 & 1.501 & 1.0886 & $\exp(10^9)$  \\
45.18 & 1.501 & 1.0903 & $\exp(10^{10})$  \\
\hline
\hline
\end{tabular}
\end{center}
\end{table}

\begin{table}[!htb]
\caption{More asymptotically stringent bounds on the first Chebyshev function $\vartheta(x)$ of the de la Valle Poussin form valid on restricted regions~\cite{dlVP}.
(Based on reference~\cite{Johnston-Yang}.) }\smallskip
\begin{center}
\begin{tabular}{||c|c|c|c||}
\hline
\hline
$a$ & $b$ & $c$ & $x_0$  \\
\hline
\hline
357 & 0 & 1/2  & $\exp(3000)$\\
320 &  0 &1/2  & $\exp(4000)$  \\
295 &  0 &1/2  & $\exp(5000)$  \\
274 &  0 &1/2  & $\exp(6000)$  \\
263 &  0 &1/2 &  $\exp(7000)$  \\
252 &  0 &1/2  & $\exp(8000)$  \\
243 &  0 &1/2 & $\exp(9000)$  \\
249 &  0 &1/2  & $\exp(10000)$  \\
\hline
\hline
644 &  0 &1/2  & $\exp(10^5)$  \\
348 &  0 &1/2 & $\exp(10^6)$  \\
312 &  0 &1/2 & $\exp(10^7)$  \\
301 &  0 &1/2 &  $\exp(10^8)$  \\
298 &  0 &1/2 &  $\exp(10^9)$  \\
297 &  0 &1/2 & $\exp(10^{10})$  \\
\hline
\hline
1642333 &  0 &1 & $\exp(10^6)$  \\
165152 &  0 &1 & $\exp(10^7)$  \\
101831 &  0 &1 &  $\exp(10^8)$  \\
87551 &  0 &1 &  $\exp(10^9)$  \\
83063 &  0 &1 & $\exp(10^{10})$  \\
\hline
\hline
\end{tabular}
\end{center}
\label{default}
\end{table}%

\clearpage
Herein we shall show how to convert these effective bounds on $\vartheta(x)$ into effective bounds on the prime gaps $g_n = p_{n+1}-p_n$.  Specifically we shall establish both
\begin{equation}
{g_n\over p_n}  < { 2a  \;(\ln p_n)^b \; \exp\left(-c\; \sqrt{\ln p_n}\right) \over
1- a  \;(\ln p_n)^b \; \exp\left(-c\; \sqrt{\ln p_n}\right)};
\qquad (x \geq x_*);
\end{equation}
and
\begin{equation}
{g_n\over p_n}  <  3a  \;(\ln p_n)^b \; \exp\left(-c\; \sqrt{\ln p_n}\right) ;
\qquad (x \geq x_*);
\end{equation}
for some effective computable $x_*$. 
In all cases it is the presence of the exponential factor that is central to making these bounds interesting and relatively stringent.

\section{Strategy}

Let us now develop some effective bounds on prime gaps $g_n = p_{n+1}-p_n$, starting from effective bounds on the first Chebyshev function of the form
\begin{equation}
\label{E:key}
|\theta(x)-x| <  a \;x \;(\ln x)^b \; \exp\left(-c\; \sqrt{\ln x}\right); 
\qquad (x \geq x_0).
\end{equation}
For convenience rewrite our bound on the first Chebyshev function in the form
\begin{equation}
|\theta(x)-x| <  x \; f(x); 
\qquad (x \geq x_0).
\end{equation}
Here $f(x) = a (\ln x)^b \exp(-c\sqrt{\ln x})$ is easily verified to be monotone decreasing for $x > x_{peak} =\exp([2b/c]^2)$, where it takes on the value $f_{peak} = a (2b/c)^{2b} \exp(-2b)$.

Define
\begin{equation}
x_* = \max\left\{ x_0, \; \exp\left(\left[2b\over c\right]^2\right) \right\}.
\end{equation}
Then in the range $x\geq x_*$ the inequality (\ref{E:key}) is valid with $f'(x)\leq 0$.
This will be the primary range of interest for the following computations. Note that in the limit $b\to0$, appropriate to effective bounds of the de la Valle Poussin form, one has
\begin{equation}
x_* \to \max\left\{ x_0, \; 1 \; 
\right\} = x_0.
\end{equation}

Let us now take any $\epsilon\in(0,1)$ and consider the inequality
\begin{equation}
\theta(p_{n+1}-\epsilon) - (p_{n+1}-\epsilon) > 
-  (p_{n+1}-\epsilon )\; f (p_{n+1}-\epsilon );
\end{equation}
Thence
\begin{equation}
p_{n+1} < \theta(p_{n})+\epsilon + (p_{n+1}-\epsilon )\;f (p_{n+1}-\epsilon ).\end{equation}
But since this holds for all $\epsilon\in(0,1)$ we can in particular 
consider the limit $\epsilon\to 0$ and so deduce
\begin{equation}
p_{n+1} \leq \theta(p_{n}) +  p_{n+1}\; f(p_{n+1}).
\end{equation}
On the other hand from
\begin{equation}
\theta(p_{n}) - p_{n} <
p_{n}\; f(p_{n});
\end{equation}
we deduce
\begin{equation}
 p_{n} >  \theta(p_{n}) 
- p_{n}\; f(p_{n});
\end{equation}
Thence we can bound the prime gaps as
\begin{equation}
g_n < 
p_{n+1}\; f(p_{n+1})
+
p_{n}\; f(p_{n}).
\end{equation}
We now have two options:
\begin{itemize}
\item 
If $f_{peak}\leq 1$ then use $p_{n+1}= p_n + g_n$,  and the fact that $f(x)$ is monotone decreasing in the range of interest, to deduce
\begin{equation}
g_n <   (2 p_n+ g_n) \;f(p_{n}),
\end{equation}
Rearranging, and using the fact that $f(x)<1$ in the range of interest, we see
\begin{equation}
{g_n\over p_n}  <  {2 \; f(p_{n})\over 1 - \;f(p_{n})}; 
\qquad(f_{peak}\leq 1).
\end{equation}
\item If $f_{peak}>1$ it is more useful to use the standard Bertrand--Chebyshev theorem $p_{n+1} < 2 p_n$, and the fact that $f(x)$ is monotone decreasing in the range of interest, to deduce
\begin{equation}
{g_n\over p_n} <   3 \;f(p_{n}); \qquad(f_{peak}\hbox{ arbitrary}).
\end{equation}
\end{itemize}
We can summarize this in a simple Lemma.

\paragraph{Lemma:} 
\emph{
If one has somehow established a bound of the form
\begin{equation}
|\theta(x)-x| <  a \; x \; (\ln x)^b \; \exp\left(- c\; \sqrt{\ln x}\right); 
\qquad (x \geq x_0);
\end{equation}
as in Tables I, II, and III above, then 
defining 
\begin{equation}
x_* = \max\left\{ x_0, \; \exp\left(\left[2b\over c\right]^2\right) \right\},
\end{equation}
for the prime gap 
$g_n = p_{n+1}-p_n$ one has the bounds
\begin{equation}
{g_n\over p_n}  < { 2a  \;(\ln p_n)^b \; \exp\left(-c\; \sqrt{\ln p_n}\right) \over
1- a  \;(\ln p_n)^b \; \exp\left(-c\; \sqrt{\ln p_n}\right)};
\qquad (x \geq x_*; \; f_{peak}\leq 1);
\end{equation}
\begin{equation}
{g_n\over p_n}  <  3 a  \;(\ln p_n)^b \; \exp\left(-c\; \sqrt{\ln p_n}\right) ;
\qquad (x \geq x_*; \; f_{peak} \hbox{ arbitrary});
\end{equation}
These bounds certainly hold for $x \geq x_*$, but if $x_*$ is sufficiently small one might be able to widen the range of applicability to some $x \geq x_{**}$, with $x_{**} \leq x_*$, by explicit computation. 
}
\section{Effective bounds on the prime gaps}

\subsection{Widely applicable bounds}

For some widely applicable bounds of the form
\begin{equation}
{g_n\over p_n}  < { 2a  \;(\ln p_n)^b \; \exp\left(-c\; \sqrt{\ln p_n}\right) \over
1- a  \;(\ln p_n)^b \; \exp\left(-c\; \sqrt{\ln p_n}\right)};
\qquad (p_n \geq x_{**}; \; f_{peak}\leq 1);
\end{equation}
consider Table IV below. For any collection of coefficients $\{a,b,c\}$ one first calculates $x_{peak}$ and checks that $f_{peak}<1$. From that and $x_0$ one determines $x_*$. Finally, for $x_*$ sufficiently small, one determines $x_{**}$ by explicit computation.

\begin{table}[!h]
\caption{Some widely applicable effective bounds 
on the relative prime gap $g_n/p_n$.\\
Compare with parts of Table I.}\smallskip
\begin{center}
\begin{tabular}{||c|c|c|c||c|c|c|c||}
\hline
\hline
\hline
$a$ & $b$ & $c$ & $x_0$ & $x_{peak}$ & $f_{peak}$ & $x_*$ & $x_{**}$ 
\\
\hline
\hline
\hline
 0.2196138920& 1/4 & 0.3219796502 & 101 &11.15042039 &0.1659905476&101 & 11 
 \\
\hline
1/4& 1/4 & 1/4 & 31 & 54.59815003&0.2144409711 & 55 & 11
 \\
\hline
\hline
0.2428127763 & 1/4 &0.3935970880 &149 &5.021606990 & 0.1659905476&  149 & 11 \\
\hline
1/4& 1/4 & 1/3 & 43 &9.487735836 & 0.1857113288 & 43 & 11
 \\
\hline
\hline
\hline
 0.3510691792&0&  1/4  & 101 & 1 & 0.3510691792 & 101 & 2
\\
\hline
\hline
0.2748124978 &0&  1/4 &  149 & 1 & 0.2748124978  & 149 & 11
\\
\hline
0.4242102935 &0&  1/3 &  149 & 1 &0.4242102935 & 149 & 2
\\
\hline
\hline
\hline
1 & 0 & 1/4 & 2 & 1 & 1& 2 & 2\\
\hline
1 & 0 & 1/3 & 3 & 1 & 1 & 3 & 2\\
\hline
\hline
1/2 & 0 & 1/4 & 29 &1 & 1/2 & 29 & 2\\
\hline
1/2 & 0 & 1/3 & 41 & 1 & 1/2 & 41 & 2\\
\hline
\hline
\end{tabular}
\end{center}
\end{table}

\subsection{Some intermediate strength bounds}

Now consider some intermediate strength bounds, (now trading off the range of applicability versus tightness of the bound), based on the Fiori--Kadiri--Swidinsky~\cite{Fiori-et-al} and Johnston--Yang~\cite{Johnston-Yang} results.  Consider the coefficients presented in Table V, applied to bounds  of the form
\begin{equation}
{g_n\over p_n}  <  3 a  \;(\ln p_n)^b \; \exp\left(-c\; \sqrt{\ln p_n}\right) ;
\qquad (x \geq x_*; \;\; f_{peak} \hbox{ arbitrary});
\end{equation}
For any collection of coefficients $\{a,b,c\}$ one first calculates $x_{peak}$, (and also verifies $f_{peak}>1$). From that and $x_0$ one determines $x_*$, which is sometimes distressingly large. Finally one determines $x_{**}$ by direct computation. Unfortunately the resulting bounds, while widely applicable, are not particularly stringent. 

\begin{table}[!h]
\caption{Some intermediate strength effective bounds 
on the relative prime gap $g_n/p_n$. \\Based on Fiori--Kadiri--Swidinsky~\cite{Fiori-et-al} and Johnston--Yang~\cite{Johnston-Yang}.\\
Compare with parts of Table I.}\smallskip
\begin{center}
\begin{tabular}{||c|c|c|c||c|c|c|c||}
\hline
\hline
\hline
$a$ & $b$ & $c$ & $x_0$ & $x_{peak}$ & $f_{peak}$ & $x_*$ & $x_{**}$ 
\\
\hline
\hline
9.220226 & 3/2 & 0.8476836 & 2 & 275108.1632 &20.34794437 & 275109 & 2\\
\hline
9.40 & 1.515 & 0.8274 & 2 & 667160.3762 &23.19042582 & 667161 & 2\\
\hline
\hline 
295 &0&  1/2 & 2 &  1 & 295 & 2 & 2
\\
\hline
\hline
385 &0&  1/2  & 2 & 1 & 385 & 2 & 2
\\
\hline
\hline
\hline
\end{tabular}
\end{center}
\end{table}
%


\subsection{Asymptotically stringent bounds}

Finally, based on Tables II and III,  consider asymptotically stringent bounds of the form
\begin{equation}
\label{E:fff}
{g_n\over p_n}  <  3 a  \;(\ln p_n)^b \; \exp\left(-c\; \sqrt{\ln p_n}\right) ;
\qquad (x \geq x_*; \; f_{peak} \hbox{ arbitrary});
\end{equation}
For any collection of coefficients $\{a,b,c\}$ one first calculates $x_{peak}$. From that and $x_0$ one determines $x_*$.
\begin{itemize}
\item For all entries in Table II it is easy to verify that $x_{peak} =\exp([2b/c]^2) \ll x_0$, (and for that matter, $f_{peak}>1$).  Thence for all entries in Table II one has $x_* = x_0$. Since $x_*$ is truly enormous direct computation of $x_{**}$ is hopeless.
In short, the effective bounds on $\vartheta(x)$ given in terms of the parameters $\{a,b,c,x_0\}$ of Table II directly imply effective bounds (\ref{E:fff}) on $g_n/p_n$ in terms of the same parameters $\{a,b,c,x_0\}$.
\clearpage
\item 
For all entries in Table III, since they are all of de la Valle Poussin form, (that is, $b=0$),  it is trivial  to verify that $x_{peak} =\exp([2b/c]^2) =1$, (and for that matter, $f_{peak}=a >1$).  Thence for all entries in Table III one trivially has $x_* = x_0$. Since $x_*$ is truly enormous direct computation of $x_{**}$ is hopeless.
In short, the effective bounds on $\vartheta(x)$ given in terms of the parameters $\{a,b,c,x_0\}$ of Table~III directly imply effective bounds (\ref{E:fff}) on $g_n/p_n$ in terms of the same parameters $\{a,b,c,x_0\}$.
\end{itemize}

\section{Conclusions}\label{S:discussion}

We have developed a number of effective bounds on the prime gaps $g_n/p_n$. 
Some of these effective bounds could in principle have been deduced almost 50 years ago. Others rely on recent numerical work from the previous decade. 
In the interests of clarity, let me quote a few explicit examples:
\begin{equation}
{g_n\over p_n} < { {1\over2} (\ln p_n)^{1/4} \exp(-\sqrt{\ln p_n}/3)
\over
1 -{1\over4} (\ln p_n)^{1/4} \exp(-\sqrt{\ln p_n}/3)};
\qquad (p_n \geq 2);
\end{equation}
\begin{equation}
{g_n\over p_n} < { \exp(-\sqrt{\ln p_n}/3)
\over
1 -{1\over2} \exp(-\sqrt{\ln p_n}/3)};
\qquad (p_n \geq 2);
\end{equation}
\begin{equation}
{g_n\over p_n} <  885 \exp(-\sqrt{\ln p_n}/2);
\qquad (p_n \geq 2);
\end{equation}
and the asymptotically tighter result
\begin{equation}
{g_n\over p_n} <  4926999 \exp(-\sqrt{\ln p_n});
\qquad (p_n \geq \exp(10^6)).
\end{equation}
In all cases it is the presence of the exponential factor that is central to making these bounds interesting and relatively stringent. 

\section*{Acknowledgements}

MV was directly supported by the Marsden Fund, \emph{via} a grant administered by the
Royal Society of New Zealand.


\clearpage


\begin{thebibliography}{99}
\newcommand{\arXiv}[1]{arXiv:~{\href{https://arxiv.org/abs/#1}{\color{blue}#1}}}


\bibitem{Schoenfeld}
Lowell Schoenfeld,\\
``Sharper bounds
for the Chebyshev functions $\theta(x)$ and $\psi(x)$. II'', \\
Mathematics of Computation, 
{\bf 30 \#134} (April 1976) 337--360.\\
\doi{10.1090/S0025-5718-1976-0457374-X}

\bibitem{Trudgian}
Tim Trudgian,
``Updating the error term in the prime number theorem'',\\
The Ramanujan Journal {\bf 39} (2016) 225--236,
\doi{10.1007/S11139-014-9656-6}.
\arXiv{1401.2689} [math.NT]


\bibitem{Johnston-Yang}
Daniel R. Johnston, Andrew Yang,\\
``Some explicit estimates for the error term in the prime number theorem'',\\
\arXiv{2204.01980} [math.NT]

\bibitem{Fiori-et-al}
Andrew Fiori, Habiba Kadiri, Joshua Swindisky,\\
``Sharper bounds for the error term in the Prime Number Theorem'',\\
\arXiv{2206.12557} [math.NT]

\bibitem{Broadbent-et-al}
S. Broadbent, H. Kadiri, A. Lumley, N. Ng, and K. Wilk. \\
``Sharper bounds for the Chebyshev function $\theta(x)$''. \\
Math. Comp. 90.331 (2021), pp. 2281--2315.
[\arXiv{2002.11068} [math.NT]]

\bibitem{Poussin}
Charles Jean de la Vall\'e Poussin,\\
``Recherches analytiques sur la th\'eorie des nombres premiers'',\\ Ann. Soc. Scient.
Bruxelles, deuxi\'eme partie {\bf20}, (1896), pp. 183--256

\bibitem{dlVP}
Matt Visser,\\
``Effective de la Vall\'e Poussin style bounds on the first Chebyshev function'', \\
\arXiv{2211.00840} [math.NT]


\bigskip
\hrule\hrule\hrule








\end{thebibliography}
\end{document}